\documentclass[12pt,a4paper]{article}
\usepackage{graphicx}
\usepackage{graphics}
\usepackage{amsmath, amssymb}
\usepackage[usenames]{color}
\usepackage{colortbl}

\newcommand{\vp}{{\mathbf p}}
\newcommand{\vu}{{\mathbf u}}
\newcommand{\vz}{{\mathbf z}}
\newcommand{\vf}{{\mathbf f}}

\newcommand{\vy}{{\mathbf y}}

\begin{document}

\sloppy

%\journal{Applied Mathematics and Computation}

%\begin{document}



\title{Ergodicity Bounds for the Markovian Queue With Time-Varying Transition Intensities, Batch Arrivals and One Queue Skipping Policy}

\author{{A. I. Zeifman\footnote{{Vologda State University;  Institute of Informatics Problems of the FRC CSC RAS; Moscow Center for Fundamental and Applied Mathematics; Vologda Research Center RAS; e-mail a$\_$zeifman@mail.ru}}},
{R. V. Razumchik\footnote{Institute of Informatics Problems of the FRC CSC RAS; Moscow Center for Fundamental and Applied Mathematics; e-mail rrazumchik@ipiran.ru}},
{Y. A. Satin\footnote{Vologda State University; e-mail yacovi@mail.ru}},
{I. A. Kovalev\footnote{Vologda State University; e-mail kovalev.iv96@yandex.ru}}}

\date{}

\maketitle

{\bf Abstract.}
In this paper we revisit the Markovian queueing system with
a single server, infinite capacity queue
and one special queue skipping policy.
Customers arrive in batches but are served one by one
in any order. The size of the arriving batch
becomes known upon its arrival and at any time instant the
total number of customers in the system is also known.
According to the adopted queue skipping policy
if a batch, which size is greater than the current total number of customers in the system,
arrives, all  customers currently residing in the system are removed from it
and the new batch is placed in the queue. Otherwise the new batch is lost
and does not have effect on the system.
The distribution of the total number of customers in the system
is under consideration under assumption that
the arrival intensity $\lambda(t)$ and/or the service intensity $\mu(t)$
are non-random functions of time.
We provide the method for the computation of the upper bounds
for the rate of convergence of system size to the limiting regime,
whenever it exists, for any bounded $\lambda(t)$ and $\mu(t)$
(not necessarily periodic) and any distribution of the batch size.
For periodic intensities $\lambda(t)$ and/or $\mu(t)$ and
light-tailed distribution of the batch size
it is shown how the obtained bounds can be used to numerically
compute the limiting distribution of the queue size
with the given error. Illustrating numerical examples are provided.

{\bf keywords:}
time-varying queueing system,  queue skipping policy,
rate of convergence bounds.

%\MSC[2010] 60J271 \sep 60J22 \sep 60J28

% \linenumbers

\section{Introduction}
The two most common viewpoints at a queueing system performance are
the point of view of the system's owner and of the system's clients. Usually their goals
are conflicting. While a client aims at minimizing one (or more)
characteristic of the jobs\footnote{Throughout this paper, when talking about a performance characteristic, we mean
its long-run value i.e. its value when the system is in the stationary or limiting regime.},
which he submits to the system (for example, job's mean response time),
the system's owner seeks to maximize the utilization of the resources (for example, processor utilization).
Both viewpoints have received attention from the operation research community
in the last decades. But performance evaluation of queueing systems
from the client's perspective seems to have attracted more attention.

\begin{figure}[h]
 \centering
  \includegraphics[width=0.8\linewidth]{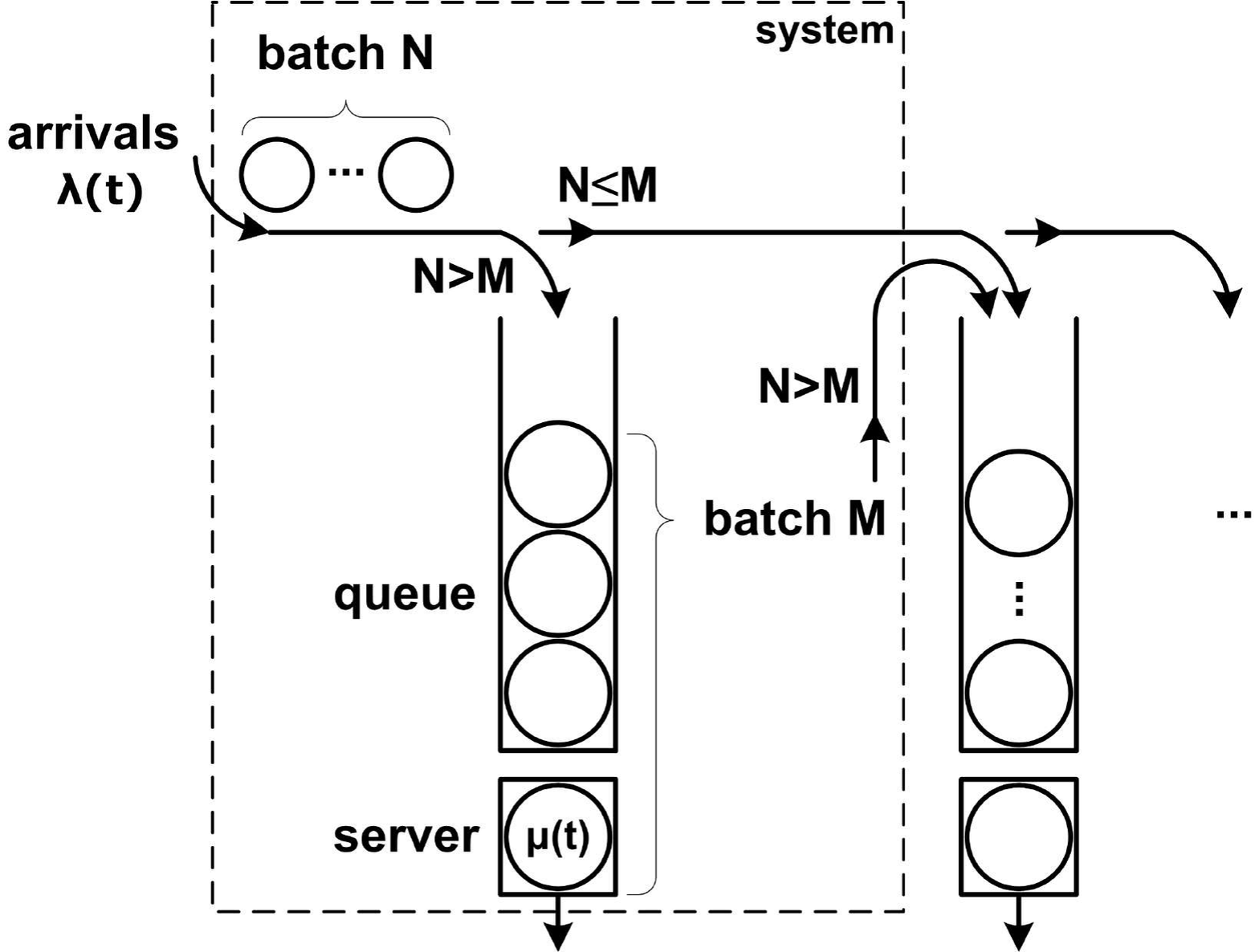}
%  \vspace{-3cm}
\caption{Model of the system with the queue skipping policy. It can be seen that the batches discarded from the system are (supposed not to be cleared but) offloaded to the next system with the same queue skipping policy and so on. This figure is the refinement of the figure~1 in \cite{elsayed}.}
\label{fig:sysplan}
\end{figure}

If we limit ourselves only to single-server queues, then
probably one of the best-known results here is the optimality of the SRPT
(shortest remaining processing time) policy with respect to
the job's (or customer's) mean response time\footnote{See \cite{schrage} and \cite{grosof}
for the latest results for multi-server queues.}. As is known, under the SRPT
at all times the server is working on the ``shortest'' job. In \cite{marin2020}
(based on the previous works \cite{balsamo,pittel}) it was noticed
that somewhat similar idea can be used to construct policies\footnote{Of course,
the requirement of the minimality of job's mean response time
under such a policy is dropped.},
which increase the utilization of all the servers in a system.
One such policy, further referred to as the ''queue skipping'' policy,
works as follows (see Fig.~\ref{fig:sysplan}).
Assume that customers arrive to the system in batches,
but are served one by one in any order. The size of the arriving batch
becomes known upon its arrival and at any time the current system size (i.e.
total number of customers in the system) is also known.
According to the queue skipping policy\footnote{As mentioned above this policy is beneficial from the viewpoint of the system's owner since, when applied to the systems in series, as in Fig.~\ref{fig:sysplan}, it increases servers' utilizations. It can also be seen from Fig.~\ref{fig:sysplan} that systems in series with such a policy are
in some sense similar to ordered-entry queues, which are well-known models for conveyor systems (open and closed-loop) with multiple unloading stations (see \cite{razumchik,disney,matsui}). } if a batch, which size is greater than the current system size,
arrives to the system, all current customers in the system are removed from it
and the new batch is placed in the queue. Otherwise the new batch is lost and does not have 
affect on the system.
In \cite{marin2020} the authors applied the time-reversed chains techniques (developed in \cite{harrison,harrison2014,kelly,marin2017}),
to study the performance of the $M/M/1$ system with such a queue skipping policy
and generally distributed batch size, and, among other results,
demonstrated the effect of the policy on the system utilization.

In this paper the effort is made to evaluate the performance of the
same system, but with time-varying intensities i.e.
when the arrival intensity $\lambda(t)$ and/or the service intensity $\mu(t)$
are non-random functions of time. Since the system is Markovian,
the (time-varying) probability density function
of the total number $X(t)$ of customers in the system 
evolves according to the system of
ordinary differential equations (ODEs) -- Kolmogorov forward equations.
Except for very special cases\footnote{For example, when the arriving
batch is always of size $1$.}, this system cannot be solved.
Moreover if the batch size distribution has infinite support
there are infinitely many ODEs in the system
and the exact analytic solution is not possible\footnote{It must be mentioned
that for practical purposes numerical computation of the
time-dependent (and limiting) densities is always possible. Indeed
$X(t)$ is the inhomogeneous continuous time birth-and-death Markov chain.
Thus whenever its state space is finite (or is somehow truncated to become finite) one can apply various uniformization algorithms (see, for example, \cite{burak}).}.
Here we adopt the approximation approach\footnote{For
probably the latest review of other approaches for time-varying queues
one case refer to \cite[Section~1]{whitt} and \cite{schwarz}.}, which circumvents
the difficulties by truncating the system of ODEs.

The contributions of this paper can be summarized as follows:
\begin{itemize}

\item We provide the method for the computation of the upper bounds
for the rate of convergence of $X(t)$ to the limiting regime,
whenever it exists, for any (not necessarily periodic)
arrival intensity $\lambda(t)$ bounded from above by a constant,
any locally integrable on $[0,\infty)$ service intensity $\mu(t)$,
and any distribution of the batch size.
This method uses the notion of the logarithmic
norm of the linear operator and is based on the previous research
\cite{granovsky,zeifman2006};

\item For periodic intensities $\lambda(t)$ and/or $\mu(t)$ and
light-tailed distribution of the batch size (which includes the geometric
distribution and distributions with finite support)
we show how one can compute the truncation
threshold $t^*$, such that the probability distribution
of $X(t)$ for $t>t^*$ ``almost forgets'' the distribution of $X(0)$.
The latter means the all performance characteristics, which depend only
on $X(t)$, can be considered as limiting characteristics for $t>t^*$
containing only a small error, which can be computed.
Since under periodic intensities the solution of the
system of ODEs governing the behaviour of $X(t)$
is also periodic, it is sufficient to
compute\footnote{If the batch size distribution has infinite support
we still have infinitely many ODEs and thus we have to perform another truncation
of the system, which introduces additional error to the final result.}
numerically the solution only in the interval $[t^*,t^*+T]$, where
$T$ is chosen manually such that the interval $[t^*,t^*+T]$ includes at
least one period of the solution;

\item Using the developed approximation approach, we numerically
compare the utilization of the system with
the queue skipping policy, periodic arrival intensity $\lambda(t)$, fixed service intensity $\mu$
and geometrically distributed batch size with mean $b$ with that of the classical
$M_t/M/1/0$ queue with the same arrival intensity $\lambda(t)$
and service intensity $\mu/b$.
This is intended to demonstrate the effect of the queue skipping policy
on the system utilization.

\end{itemize}

\noindent
Even though the obtained rate of convergence bounds do hold
for any bounded arrival intensity $\lambda(t)$,
locally integrable service intensity $\mu(t)$
and any distribution of the batch size, the proposed
method for finding the truncation threshold $t^*$ has
so far limited applicability.
This is due to the fact that for long-tailed batch size
distributions (i.e. those which have tails heavier than
the geometric distribution) so far we were unable to find the condition,
which guarantees the existence of the limiting regime of $X(t)$ (even for periodic intensities).

The paper is structured as follows. In Section 2 we repeat the description of the model.
Section~3 contains the main result (see the Theorem and Corollaries 1--3).
We demonstrate the method based on the logarithmic norm
to bound (from above) the rate of convergence of $X(t)$ to the limiting regime
(assuming that it exists). Here we also show that for
geometrically distributed batch size the limiting regime always exists.
In Section~4 the numerical example is given. Section 5 concludes the paper.

\section{System description}

Consider the $M_t/M_t/1$ queue with intensities being periodic functions of time
and the queue skipping policy.
Customers arrive to the system in batches according
to the inhomogeneous Poisson process with intensity $\lambda(t)$.
The size of an arriving batch becomes known upon its arrival 
and is the random variable with the given probability distribution $\{ b_n, n \ge 1 \}$,
having finite mean $\sum_{k=1}^\infty B_k$, $B_k = \sum_{n=k}^{\infty}b_n$.
The adopted queue skipping policy implies that
whenever a batch arrives to the system its size, say $\widehat B$, is compared with the
current total number of customers in the system, say $\widetilde B$.
If ${\widehat B} > { \widetilde B}$, then all customers, which are currently in the system,
are instantly removed from it, and the whole batch ${\widehat B}$ is placed in 
the queue and the first customer in the batch enters server.
If ${\widehat B} \le  { \widetilde B}$ the new batch leaves the system without having any
effect on it.
Whenever the server becomes free one customer from the queue (if there is any) enters server\footnote{Since we do not study waiting time characteristics in this paper,
the service discipline is unimportant and
for certainty one can consider that customers are served in FIFO or LIFO or RANDOM order.} and
gets served according to exponential distribution with intensity $\mu(t)$.

\section{Main result}

Let $X(t)$ be the total number of customers in the system at time $t$.
From the system description it follows that $X(t)$ is
the Markov chain with continuous time and discrete state space
$\mathcal{X}=\{ 0, 1, 2, \dots, b^*\}$,
where $b^*$ is the maximum possible batch size i.e. $b^*=\max_{n\ge 1} (b_n>0)$.
If the batch size distribution has infinite support
then the state space $\mathcal{X}$ is countable;
otherwise it is finite.

Denote by $Q(t)$ the intensity matrix (infinitesimal generator) of $X(t)$.
It is straightforward to see that $Q(t)$ has the form

{\footnotesize
\begin{equation*}
Q(t) =\begin{pmatrix}
-\lambda(t) & \lambda(t)b_1 & \lambda(t)b_2 &   \dots \\
\mu(t) & - \left ( \mu(t) + \lambda(t)B_2 \right ) & \lambda(t)b_2   &  \dots \\
0& \mu(t) & - \left ( \mu(t) + \lambda(t)B_3 \right ) &   \dots \\
0& 0& \mu(t)  &  \dots \\
\vdots& \vdots& \vdots& \ddots \\
\end{pmatrix}.
\end{equation*}}

\noindent We assume that $\lambda(t)$ and $\mu(t)$
are arbitrary non-random functions of $t$, locally integrable
on $[0,\infty)$ and that the arrival intensity is bounded by a constant i.e.
there exists $L>0$ such that $\lambda(t)\le L < \infty$ for $t\ge 0$.

Denote by $p_i(t)=\mathbf{P}(X(t)=i)$
the probability that the Markov chain $X(t)$ is in state $i$ at time $t$.
Let $\vp(t) = (p_0(t), p_1(t), \dots)^T$ be the
probability distribution vector at time $t$.
Given any proper initial condition $\vp(0)$,
the probabilistic dynamics of
the Markov chain $X\left(t\right)$
is described by the forward Kolmogorov system
of differential equations
\begin{equation}
\label{ur01}
\frac{d}{dt}\vp(t)=A(t)\vp(t),
\end{equation}
\noindent where $A(t) = Q^T(t)$  is the transposed intensity matrix.
Throughout the paper vectors are regarded as column vectors,
${\bf 0}$ denotes the vector consisting of zeros,
$I$ denotes the identity matrix and $\cdot^{T}$ denotes the matrix transpose.
The sizes of matrices will be clear from the context.
The choice of vector norms will be the $l_1$-norm, that is, $\|{\vp(t)}\|=\sum_{i\in \mathcal{X}} |p_i(t)|$; the operator norm will be the one induced by the $l_1$-norm on row vectors,
that is, $\|A(t)\| = \sup_{j \in \mathcal{X}} \sum_{i \in \mathcal{X}} |a_{ij}(t)|$.

Recall that a Markov chain $X(t)$ is called
weakly ergodic, if ${\|\vp^{*}(t)-\vp^{**}(t)\| \to 0}$ as $t \to
\infty$ for any initial conditions $\vp^{*}(0)$ and $\vp^{**}(0)$,
where $\vp^{*}(t)$ and $\vp^{**}(t)$ are the corresponding solutions
of (\ref{ur01}). The rate at which this
difference tends to zero is called the rate of convergence.
Below we present the method\footnote{This method is not new
and has already been applied to bound the rate of convergence
in other settings, for example, \cite{granovsky,zeifman2006,zeifman2018,zeifman2020}.
But the structure of the infinitesimal generator $Q(t)$ is different from
all those considered so far. This motivates the analysis
carried out below, since the opportunity to use
logarithmic norm to bound the rate of convergence heavily depends
on the structure of the infinitesimal generator.}
based on the logarithmic norm of a linear operator
function, which allows one to bound
from above this rate of convergence.

Using the normalization condition $p_0(t) = 1 - \sum_{i \ge 1, i \in \mathcal{X}} p_i(t)$
it can be checked that the system (\ref{ur01}) can be rewritten as follows:
\begin{equation}
\frac{d}{dt}{\vz}(t)= B(t){\vz}(t)+{\vf}(t), \label{2.06}
\end{equation}
\noindent where
$B(t)= \left(b_{ij}(t)\right)_{i,j=1}^{\infty}$,
$b_{ij}(t) = a_{ij}(t)- a_{i0}(t)$, and
$$
{\vf}(t)=\left( \lambda(t)b_{1}, \lambda(t)b_{2},\dots \right)^{T},
$$
$$
{\vz}(t)=\left(p_1(t), p_2(t),\dots \right)^{T},
$$
$$
B(t)=
\begingroup % keep the change local
\setlength\arraycolsep{1pt}
{\footnotesize \left(
\begin{array}{cccccccc}
- \!\left ( \mu(t) \!+\! \lambda(t) \right ) & \mu(t) \!-\!\lambda(t)b_{1}
& -\!\lambda(t)b_{1} & -\!\lambda(t)b_{1} & \cdots \\
0 & -\! \left ( \mu(t) \!+\! \lambda(t)B_2 \right ) & \mu(t)
\!-\!\lambda(t)b_{2} & -\!\lambda(t)b_{2} & \cdots  \\
\vdots & \vdots & \vdots & \vdots  & \ddots \\
\end{array} \right)\!.}\label{2.07}
\endgroup
$$

\noindent Note that the matrix $B(t)$ has no probabilistic meaning.
Let ${\vz}^{*}(t)$ and ${\vz}^{**}(t)$ be the solutions
of (\ref{2.06}) corresponding to (different) initial conditions
${\vz}^{*}(0)$ and ${\vz}^{**}(0)$.
Then for the vector
${ \vy}(t) = {\vz}^{*}(t)-{\vz}^{**}(t)=\left({y}_1(t), {y}_2(t), \dots\right)^T$,
which has coordinates of arbitrary signs, we have
\begin{equation}
\frac{d}{dt}{\vy}(t)= B(t){ \vy}(t). \label{hom1}
\end{equation}

\noindent Thus all bounds on the rate of convergence to the limiting regime of
$X(t)$ correspond to the same rate of convergence bounds of the solutions of the system
(\ref{hom1}). It is more convenient\footnote{Apparently this was firstly noticed in \cite{zeifman1989}.}
to study the rate of convergence using the transformed version $B^*(t)$ of $B(t)$ given
by $B^*(t)=TB(t)T^{-1}$, where $T$ is the upper triangular matrix
of the form
$$
T=\left(
\begin{array}{ccccccc}
1   & 1 & 1 & \cdots \\
0   & 1  & 1  &   \cdots  \\
0   & 0  & 1  &   \cdots \\
\vdots & \vdots & \vdots & \ddots \\
\end{array}
\right). \label{vspmatr}
$$

\noindent Let ${\vu}(t)=T{\vy}(t)=\left( {u}_1(t), u_2(t), \dots\right)^T$. Then
by multiplying \eqref{hom1} from the left by $T$
we obtain
\begin{equation}
\frac{d}{dt} {\vu}(t) = B^*(t) {\vu}(t), \label{hom11}
\end{equation}

\noindent
where ${\vu}(t)$ is the vector with the coordinates of arbitrary signs
and the matrix $B^*(t)$ has the following structure:
$$
B^*(t)=
\begingroup % keep the change local
\setlength\arraycolsep{1pt}
{\footnotesize \left(
\begin{array}{ccccccc}
-\! \mu(t)\! -\! \lambda(t) & \mu(t)
 & 0 & 0 &  \cdots\\
0 & - \!\mu(t) \!-\! \lambda(t)B_2  & \mu(t) & 0 & \cdots\\
0 & 0  & - \!\mu(t) \!-\! \lambda(t)B_3  & \mu(t) & \cdots\\
\vdots & \vdots & \vdots & \vdots &\ddots   \\
\end{array}
\right)}.\label{B^*}
\endgroup
$$

\noindent
The matrix $B^*(t)$ is essentially non-negative, i.e. all its
off-diagonal elements are non-negative for any $t \ge 0$. From this
fact it follows that,
if the initial condition $\vu(s)$ is non-negative,
then any solution $\vu(t)$ of (\ref{hom11}) is non-negative for any $0 \le s \le t$.

Now everything is ready to apply the method of logarithmic norm.
Recall that the logarithmic norm $\gamma({B}(t))$
of the operator function $B(t)$ is defined as
$$
\gamma({B}(t)) = \lim_{h \to
+0}h^{-1}\left(\|I+hB(t)\|-1\right).
$$

\noindent Denote by $V(t, s)= V(t)V^{-1}(s)$ the Cauchy operator of the equation (\ref{hom1}).
Then $\|V(t, s)\| \le e^{\, \int_s^{t} \gamma(B(u))\, du}$. For an operator function,
which maps $l_1$-vectors into $l_1$-vectors\footnote{I.e. vectors with $l_1$-norm.} and $B(t)$ is such an operator, $\gamma({B}(t))$ is expressed as:
$$
\gamma({B}(t)) = \sup_{j \in \mathcal{X}} \left(b_{jj}(t)+\sum_{i \in \mathcal{X}, i \neq j}
|b_{ij}(t)|\right). \label{lognorm1}
$$

\noindent If the matrix $B(t)$ is essentially
non-negative then
$\gamma({B}(t)) = \sup_{j \in \mathcal{X}} \left (\sum_{i \in \mathcal{X}}
b_{ij}(t) \right )$.

Let ${\{d_i, \ i \ge 1\}}$ be a sequence of positive numbers such
that $1 = d_1 \le d_2 \le \dots$ and let $D=diag(d_1,d_2,\dots)$ be
the diagonal matrix.
By putting ${\bf w}(t)=D {\bf u}(t)$ in (\ref{hom11}),  we obtain the
following equation\footnote{ Just as the matrix $B(t)$, the matrix $B^{**}(t)$ and thus the vector ${\bf w}(t)$ have no probabilistic meaning. These transformations are needed to change the structure of the initial intensity matrix $Q(t)$ and bring it to such form, which allows exact analysis of the ergodicity bounds.}
$$
\frac{d }{dt} {\bf w}(t)= B^{**}(t){\bf w}(t), \label{hom111}
$$
\noindent where $B^{**}(t)=D
B(t)^{*}D^{-1}=\left(b^{**}_{ij}(t)\right)_{i,j=1}^\infty$. Note
that $B^{**}(t)$ is also nonnegative for any $t \ge 0$.
Put
$$
\alpha_i\left(t\right)= -\sum_{j=1}^\infty b_{ji}^{**}(t),
 \ i \ge 1, \label{posit02}
$$
\noindent and let $\alpha\left(t\right)=\inf_{i \ge 1}\alpha_i\left(t\right)$.
We have that $\gamma({B^{**}}(t)) = -\alpha\left(t\right)$ and
\begin{equation*}
\|{\bf w}(t)\| \le e^{-\int_s^t \alpha\left(u\right)\,du}\|{\bf w}(s)\|,\label{lognorm12}
\end{equation*}
\noindent for any $s,t$ such that $0 \le s \le t$.

Let $\delta <1$ be a positive number, and $d_{k+1}=\delta^{-k}$, $k \ge 1$.
Then the values of $\alpha_i$ are equal to:
$$
\alpha_1\left(t\right)= \lambda(t) + \mu(t),
$$
$$
\alpha_k\left(t\right)= \lambda(t)B_k + \mu(t)\left(1-\delta\right),
 \quad k \ge 2,
$$
\noindent and therefore, we can
bound $\alpha\left(t\right)$ by the following way:
\begin{equation}
\alpha\left(t\right) \ge \alpha^*\left(t\right)=
\left(1-\delta\right) \mu(t). \label{al01}
\end{equation}

So far we have assumed that the limiting regime of $X(t)$ exists.
Its existence in the considered queue depends on the form of
the batch size distribution $\{ b_n, n \ge 1 \}$.
Below we show that when the tail of the distribution
is geometric or lighter then the limiting regime of $X(t)$
always exists, whereas for heavier tails the question remains open.
We start with the analysis of (\ref{2.06}).
Let $V(t,s)$ be the Cauchy operator for (\ref{2.06}).
Then
%\begin{equation}
$$
\vz(t) =  V(t)\vz(0) + \int_0^t {V}(t,\tau) \vf(\tau)\, d\tau.
\label{2.071}
$$
%\end{equation}

Put ${\bf r}(t)=DT\vz(t)$. Then instead of (\ref{2.06}) we get:
\begin{equation}
\frac{d}{dt}{\bf r}(t)= B^{**}(t){\bf r}(t)+{\bf f^{**}}(t),
\label{2.072}
\end{equation}
\noindent where ${\bf f^{**}}(t)=DT{\bf f}(t)$, and
\begin{equation}
{\bf r}(t) =  V^{**}(t){\bf r}(0) + \int_0^t {V}^{**}(t,\tau) {\bf
f^{**}}(\tau)\, d\tau, \label{2.073}
\end{equation}
\noindent where $\|V^{**}(t,s)\| \le e^{-\int_s^t \alpha^*\left(u\right)\,du}$
due to (\ref{al01}). If there exist $N>0$ and $a>0$ such that
\begin{equation}
e^{-\int_s^t\alpha^*\left(u\right)\,du} \le   N e^{-a(t-s)},
\label{2.075}
\end{equation}
\noindent for any $0 \le s \le t$,
then we have exponential weak ergodicity in the corresponding norm.

Let there exist $0<q<1$ and $C>0$ such that $b_k \le C q^k$ for all $k \ge 1$.
Then
\begin{eqnarray}
\|{\bf r}(t)\| \le \| V^{**}(t)\| \|{\bf r}(0)\| + \nonumber
\\ + \int_0^t
\|{V}^{**}(t,\tau)\| \|{\bf f^{**}}(\tau)\|\, d\tau \le
 N e^{-at}\|{\bf r}(0)\| + \nonumber
\\ + \int_0^t
N e^{-a(t-\tau)}K\, d\tau \le \frac{NK}{a}+ o(1), \label{2.076}
\end{eqnarray}
\noindent where
\begin{eqnarray} \|{\bf f^{**}}(t)\| = \|DT{\bf
f}(t)\|= \nonumber
\\ =C\lambda(t)\left(\frac{d_1 q}{1-q}+\frac{d_2
q^2}{1-q}+\dots\right) \le \nonumber
\\ \le CL
\left(\frac{q}{1-q}+\frac{\delta^{-1} q^2}{1-q}+\frac{\delta^{-2}
q^3}{1-q}+ \dots\right)= \nonumber
\\=  \frac{CL\delta q}{(\delta-q)(1-q)}=K,\label{2.0761}
\end{eqnarray}
\noindent
for $\delta > q$.

Assume now that $b_k \to 0$  more slowly, say $b_k \ge k^{-s}$ for some  $s >1$.
Firstly note that essential non-negativity of $B^{**}(t)$ implies
non-negativity of the matrix ${V}^{**}(t,s)$ for any $0 \le s \le
t$. Now if ${\bf r}(s) \ge {\bf 0}$ then ${\bf r}(t) \ge {\bf 0}$
for any $t \ge s$. This follows from non-negativity of ${\bf
f^{**}}(t)$ and ${V}^{**}(t,s)$ in (\ref{2.073}).

Put ${\bf r}(0) \ge {\bf 0}$. Then $\|{\bf r}(t)\| = \sum_{k \in  \mathcal{X}} r_k(t)$,
for any $t \ge 0$. From (\ref{2.072})
for any $t \ge 0$ and for any $0<\delta <1$ we obtain
\begin{eqnarray*}
\label{2.080}
\frac{d}{dt}{\|\bf r}(t)\|= \sum_{k \in \mathcal{X}}\frac{d}{dt}{ r_k}(t) \ge \|{\bf
f^{**}}(t)\| = \nonumber \\
= \lambda(t)\left(\left(b_1+b_2+b_3+\dots\right)+  \delta^{-1}\left(b_2+b_3+\dots\right)+ \right.
 \\
\left.+\delta^{-2}\left(b_3+\dots\right) +\dots \right) \ge \lambda(t)\sum_{k \ge 1} \delta^{-k} k^{-s}= \infty. \nonumber
\end{eqnarray*}
\noindent
Hence the corresponding equation (\ref{2.072}) does not have a
limiting solution in the corresponding norm.
Thus weak ergodicity and the existence of
limiting characteristics is guaranteed only if the tail of the
batch size distribution is geometric or lighter.
We summarize the findings in the following

\medskip

\textit{{\bf Theorem.} Assume that
exist $0<q<1$ and $C>0$ such that $b_k \le C q^k$ for all $k \ge 1$, and that
\begin{equation}
\int_0^{\infty} \mu(t)\,dt = + \infty. \label{lognorm21}
\end{equation}
Then the Markov chain $X(t)$ is  weakly ergodic and for any
initial condition ${\bf w}(0)$ and any  $t \ge 0$ the following
upper bound holds:
$$
\|{\bf w}\left(t\right)\| \le e^{-\int_0^t
{\left(1-\delta\right)\mu\left(u\right)\,du}}\|{\bf w}(0)\|, \label{t001}
$$
\noindent for any $\delta \in (q,1)$.
If (\ref{2.075}) holds for some $N>0$ and $a>0$,
then $X(t)$ is exponentially weakly ergodic.}

Now we can obtain the bounds in more natural norms.
Firstly note that $ \|\vp^*(t)-\vp^{**}(t)\| \le 2 \|\vz^*(t)-\vz^{**}(t)\| \le
4\|{\bf w}(t)\|$ and $\|\vz (t)\|_{1E}\le W^{-1}\|\bf w(t)\|$ (see \cite{zeifman2006}),
where $l_{1E}=$ \\ $\left\{z(t)=(p_1(t),p_2(t),\ldots)^T:
\|z(t)\|_{1E}\equiv\sum_{n \in \mathcal{X}} n |p_n(t)| < \infty\right\}$
and $W=\inf_{k\ge 1}\frac{d_{k}}{k} > 0$.

\medskip

\textit{{\bf Corollary 1.}
Under the assumptions of the Theorem the Markov chain $X(t)$
has a limiting mean, say $\phi(t)$,
and the following rate of convergence bounds hold:
\begin{equation}
\|\vp^*(t)-\vp^{**}(t)\| \le 4 e^{-\int_0^t
{\left(1-\delta\right)\mu\left(u\right)\,du}}\|{\bf w}(0)\|, \label{t002}
\end{equation}
\begin{equation}
|E(t,k) - \phi(t)| \le \frac{4}{W} e^{-\int_0^t
{\left(1-\delta\right)\mu\left(u\right)\,du}}\|{\bf w}(0)\|, \label{t003}
\end{equation}
\noindent
where $E(t,k)=\sum_{n \in \mathcal{X}} n p_n(t)$ is the mean number of customers
in the system at time $t$, given that initially there where $k$ customers
in the system i.e. $p_k(0)=1$.}

\medskip

\textit{{\bf Corollary 2.}
Let $X(t)$ be a homogeneous Markov chain i.e.
$\lambda(t)=\lambda$ and $\mu(t)=\mu$.
Then $X(t)$ is strongly ergodic and for any initial
condition ${\bf w}(0)$ and any  $t \ge 0$ the
following upper bounds hold:
%\begin{eqnarray}
$$
\|{\bf w}\left(t\right)\| \le e^{-
{\left(1-\delta\right)\mu t}}\|{\bf
w}(0)\|, \label{t011}
$$
%\end{eqnarray}
$$
\|\vp^*(t)-\vp^{**}(t)\| \le 4 e^{-
{\left(1-\delta\right)\mu t}}\|{\bf
w}(0)\|, \label{t012}
$$
$$
|E(t,k) - \phi(t)| \le \frac{4}{W} e^{-
{\left(1-\delta\right)\mu t}}\|{\bf
w}(0)\|. \label{t013}
$$
}

\textit{{\bf Corollary 3.}
Let the arrival  intensity $\lambda(t)$ and the service intensity $\mu(t)$ be
$1-$periodic. Then the assumptions of Theorem are equivalent to the
inequality $\int_0^1 \mu(t)\, dt > 0$.
Moreover the limiting probability distribution of the
Markov chain $X(t)$ is $1-$periodic and the limiting
mean is $1-$periodic as well.}

From the Theorem and Corollaries 1--3 it follows that
that the bounds on the rate of convergence hold
for common intensity functions.
If the latter are periodic in time then the
limiting probability characteristics of the $X(t)$
(whenever they exists) are also periodic.

\section{Numerical example}

In all the examples presented below it is assumed that the batch size distribution
$\{ b_k, \ k \ge 1\}$ is geometric i.e. $b_k=(1-q)q^{k-1}$, $k \ge 1$, $0<q<1$.
It is also assumed that the arrival and/or transition intensities
are periodic functions. Given that $\{ b_k, \ k \ge 1\}$ is geometric,
irrespective of the parameter of the geometric distribution,
periodic intensities guarantee the existence of the
(periodic) limiting distribution of $X(t)$.

When the intensities are periodic but  the batch size
has a general distribution we were unable so far to find
even sufficient conditions of the existence of the limiting distribution of $X(t)$.

Below we consider two examples: one is devoted to the discussion of the convergence
bounds obtained and in the other the properties of the queue skipping policy are
illustrated.

\subsection{Example 1}

In this example it is demonstrated how exactly the upper bound on the rate of
convergence of $X(t)$ can be computed. Let
both the arrival and service intensities be periodic
and equal to $\lambda(t)=1+\sin (2 \pi t)$ and $\mu(t)=1+\cos (2 \pi t)$.
 Let $q=\frac{2}{3}$, i.e. the batch size distribution be $b_k=\frac{2^{k-1}}{3^{k}}$, $k \ge 1$,
i.e. the mean batch size $\sum_{k=1}^\infty k b_k$ is $3$.
From (\ref{2.076}) it follows that in order to compute the upper
bound on the rate of convergence, one needs to
choose firstly $\delta \in (q,1)$ and secondly $\|{\bf w}(0)\|$.
Namely, on the one hand, for the better rate of convergence we should choose smallest possible  $\delta$ , and on the other hand, for better bounding of  $\|{\bf w}(0)\|$ we should choose as much $ \delta $ as possible.
Put   $\delta =\frac{5}{6}$.
Thus $\alpha^*(u)=\frac{1}{6} \mu(t)$
and from (\ref{2.075}) it follows that
$$
e^{-\int_s^t\alpha^*(u)\,du}=
e^{-\frac{1}{6}\int_s^t (1 +  \cos (2 \pi u) ) \,du}
 \le 2 e^{- \frac{1}{6} (t-s)},
$$
hence in the right part of (\ref{2.075}) we can put $a=\frac{1}{6}$ and $N=2$.
Now let us consider the choice of $\|{\bf w}(0)\|$. 

Consider  (\ref{2.0761}). We have $L=2$ and $C=\frac{1}{3}$. Thus $K =\frac{20}{3}$ in (\ref{2.0761}) and from (\ref{2.076}) it follows that $\limsup_{t \to \infty}\|{\bf r}(t)\| \le 80$.
Since ${\bf w}(t)= DT{\bf y}(t)$, the inequality (\ref{2.076})
guarantees that the $l_1$-norm of the limiting distribution of $X(t)$
does not exceed 80 i.e. $\|{\bf w}(0)\| \le 80$.
Thus (\ref{t002}) gives
\begin{eqnarray}
\|\vp^*(t)-\vp^{**}(t)\| \le 160 e^{- \frac{t}{6}}
\label{ex1-1}
\end{eqnarray}
for any initial conditions $\vp^*(0)$ and $\vp^{**}(0)$.
For example, if $t=80=t^*$ then the right part of
(\ref{ex1-1}) does not exceed $10^{-3}$
i.e. starting from $t>t^*$ the system ``forgets'' its initial state and
the probability distribution $X(t)$ for $t>t^*$ can be regarded
as the limiting distribution of $X(t)$.
The error (in $l_1$-norm), which is thus made is not greater than $10^{-3}$.
Moreover, since the limiting distribution of $X(t)$
is periodic, we are allowed to
solve the system of ODEs only in the interval $[0,t^*+1]$.
The probability distribution of $X(t)$ in the interval
$[t^*,t^*+1]$ is the estimate (with error not greater than $10^{-3}$ in $l_1$-norm)
of the limiting probability distribution of $X(t)$.
It must be noticed that since $b_k>0$ for all $k$,
the system of ODEs contains infinite number of equations.
Thus in order to solve it numerically one has to truncate the system.
We perform this truncation according to the method in \cite{zeifman2017}.

The upper bound on the rate of convergence
of the mean number of customers in the system $E(t,k)$ to its limiting value $\phi(t)$
is computed in the same manner.
Firstly recall that $d_{k+1}=\delta^{-k}$ and since $\delta=\frac{5}{6}$ has been
fixed above, then $d_{k+1}=\delta^{-k}=\left(\frac{6}{5}\right)^k$.
Thus $W=\inf_{k\ge 1}\frac{d_{k}}{k} =\frac{3}{4}$.
Now consider (\ref{t003}).
Thus $\|{\bf w}(0)\| \le 80$
and from (\ref{t003}) it follows that
\begin{eqnarray}
|E(t,k) - \phi(t)| \le  \frac{640}{3} e^{-\frac{t}{6}} \label{ex1-2}
\end{eqnarray}

\noindent
for any initial condition $X(0)=k$, $k \ge 0$.
Thus for $t>t^*$ the value of $E(t,k)$ can be regarded as
the limiting value of the mean number of customers
and contains the error (in $l_1$-norm) not greater than $10^{-3}$.

In Fig.~2 and Fig.~3 one can see the graphs of $p_0(t)$
and $E(t,0)$ in the interval $[0,t^*=80]$.
The ODEs are solved with the initial condition $X(0)=0$
i.e. the system is initially empty.
It can be seen that the obtained upper bounds
(\ref{ex1-1}) and (\ref{ex1-2}) are not tight:
the systems enters periodic limiting regime
before $t^*$.

\subsection{Example 2}

This example is devoted to the illustration of the properties of the queue skipping policy
in the case when the transition intensities
are time-dependent (purely Markov case is studied in \cite{marin2020}).
For simplicity we assume that only the arrival intensity $\lambda(t)$
depends on time and the service intensity is constant i.e.
$\mu(t)=\mu$, $t\ge 0$.
Since the limiting regime always exists when the batch size distribution is geometric,
no restrictions (additional to those required by the Theorem) are imposed on the arrival intensity
$\lambda(t)$.

In Fig. \ref{fig:06}, \ref{fig:08} and \ref{fig:10} one can see
how the limiting probabilities $p_0(t)$, $p_1(t)$, $p_2(t)$
and $p_3(t)$ behave depending on the service intensity $\mu$.
It is assumed that $\lambda(t)=0.8+0.1\sin (2 \pi t)$,
$b_k=2^{-k}$, $k \ge 1$, i.e. the mean batch size is $2$.
Due to the low amplitude in the arrival intensity $\lambda(t)$,
the amplitude of the limiting probabilities is also low
and is dependent on the service intensity.

As in many other queueing systems, the idle
probability $p_0(t)$ is one of the key performance indicators.
In Fig. \ref{fig:11} and \ref{fig:13}
one can see the behaviour of limiting value of $p_0(t)$,
when the service intensity is fixed ($\mu(t)=\mu=1$).
From the figures it can be seen that, as expected,
the idle limiting probability tends to $0$ when the
batch size or the arrival intensity increases.

Finally it is also of interest to compare the
limiting idle probability $p_0(t)$ of the considered
system with the queue skipping policy
with the limiting idle probability
in the pure blocking system i.e. $M_t/M/1/0$ queue
under the same arrival intensity $\lambda(t)$.
Since the arriving batch has mean size $(1-q)^{-1}$
we have to change the service intensity in the
blocking system to $\mu (1-q)$.
In Fig. \ref{fig:21} and \ref{fig:23} one can see
the graphs of $p_0(t)$ in these two systems
given that $\lambda(t)=0.8+0.1\sin (2 \pi t)$ and $\mu=1$.

We observe that even the time-inhomogeneous system with the queue skipping policy
(just like the homogeneous one studied in \cite{marin2020})
gives a much better utilisation than the pure blocking time-inhomogeneous system,
when the mean batch size is large (i.e. $q$ is close to 1)
and the arrival intensity is high.

\section{Conclusion}

Even though we have limited the discussion only to
the geometric batch size, from Section~3 it
can be seen that the method based on the logarithmic
norm of a limiting operator allows us to upper bound
the rate of convergence for any batch size distribution.
The only open problem, which persists,
is to find the conditions,
which guarantee the existence of the limiting
regime for any batch size distribution.
So far we have not been able to do it but we believe
that it is only the matter of the proper analytic point of view.
This hope is supported by the fact that in time-homogeneous
case such conditions are known (see \cite{marin2020}).

Although it is not mentioned above, being able to compute (approximately)
the limiting mean of $X(t)$ allows one to use the time-varying Little's law
to compute the average sojourn time in the system (before a customer leaves the queue
either due to an arrival or due to service completion).

The obtained results also show that in order to obtain the bounds on the rate of convergence
one does not need to know the exact values of $\lambda(t)$ and $\mu(t)$.
Instead it is sufficient to know the values of the integrals of type
(\ref{lognorm21}) i.e. the time-average intensities
$\overline{\lambda} = \frac{1}{t} \lim_{t \to \infty}\int_0^t \lambda(u)du$
and
$\overline{\mu} = \frac{1}{t} \lim_{t \to \infty}\int_0^t \mu(u)du$
(and not the exact values of $\lambda(t)$ and $\mu(t)$).
In case of 1-periodic intensities $\lambda(t)$ and $\mu(t)$
the values $\overline{\lambda}$ and $\overline{\mu}$ are exactly
the average arrival and service intensity over one period.

\bigskip

{\bf Acknowledgement.} This research was supported by Russian Science Foundation under grant 19-11-00020.

\renewcommand{\refname}{References}

\newpage
%%%%%%%%%%%%%%% -1

%%%%%%%%%%%%%%% 0

\begin{figure}[h]
 \centering
  \includegraphics[width=0.5\linewidth]{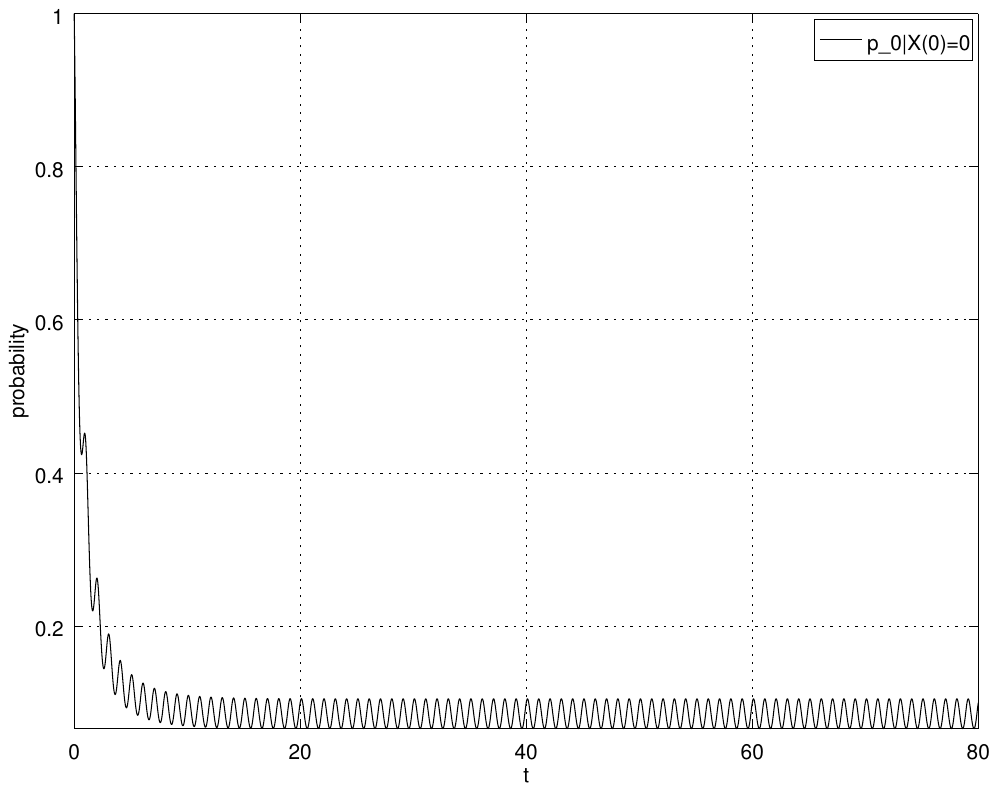}
%  \vspace{-3cm}
\caption{Rate of convergence of the empty system probability $p_0(t)$ in the interval $[0,80]$.}
%\Description{The plot for Example 1}
\label{fig:01}       % Give a unique label
\end{figure}

\begin{figure}[h]
\centering
  \includegraphics[width=0.5\linewidth]{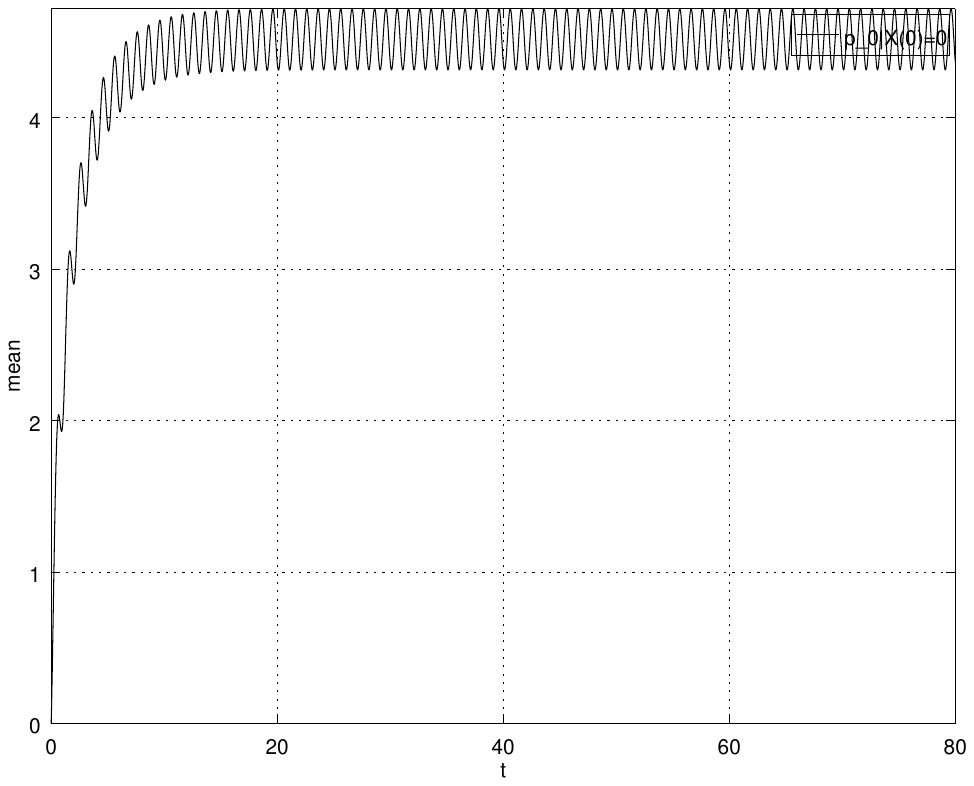}
%  \vspace{-3cm}
\caption{Rate of convergence of the mean
number of customers $E(t,0)$ in the system in the interval $[0,80]$.}
%\Description{The plot for Example 1}
\label{fig:02}       % Give a unique label
\end{figure}

%%%%%%%%%%%%%%% 1

\begin{figure}[h]
\centering
  \includegraphics[width=0.5\linewidth]{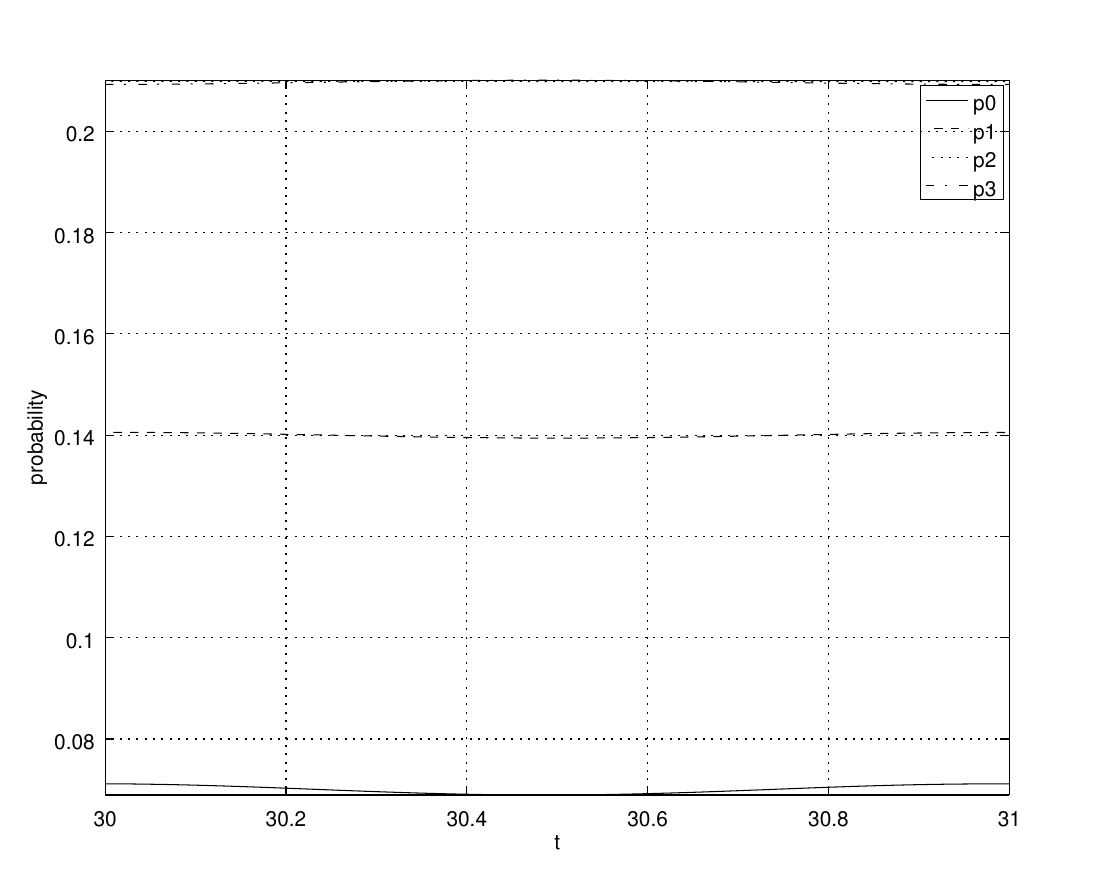}
%\vspace{-3cm}
\caption{The limiting probabilities $p_i(t)$, $0 \le i \le 3$, for $\mu =0.4$.}
%\Description{The plot for Example 2}
\label{fig:06}       % Give a unique label
\end{figure}

\begin{figure}[h]
\centering
  \includegraphics[width=0.5\linewidth]{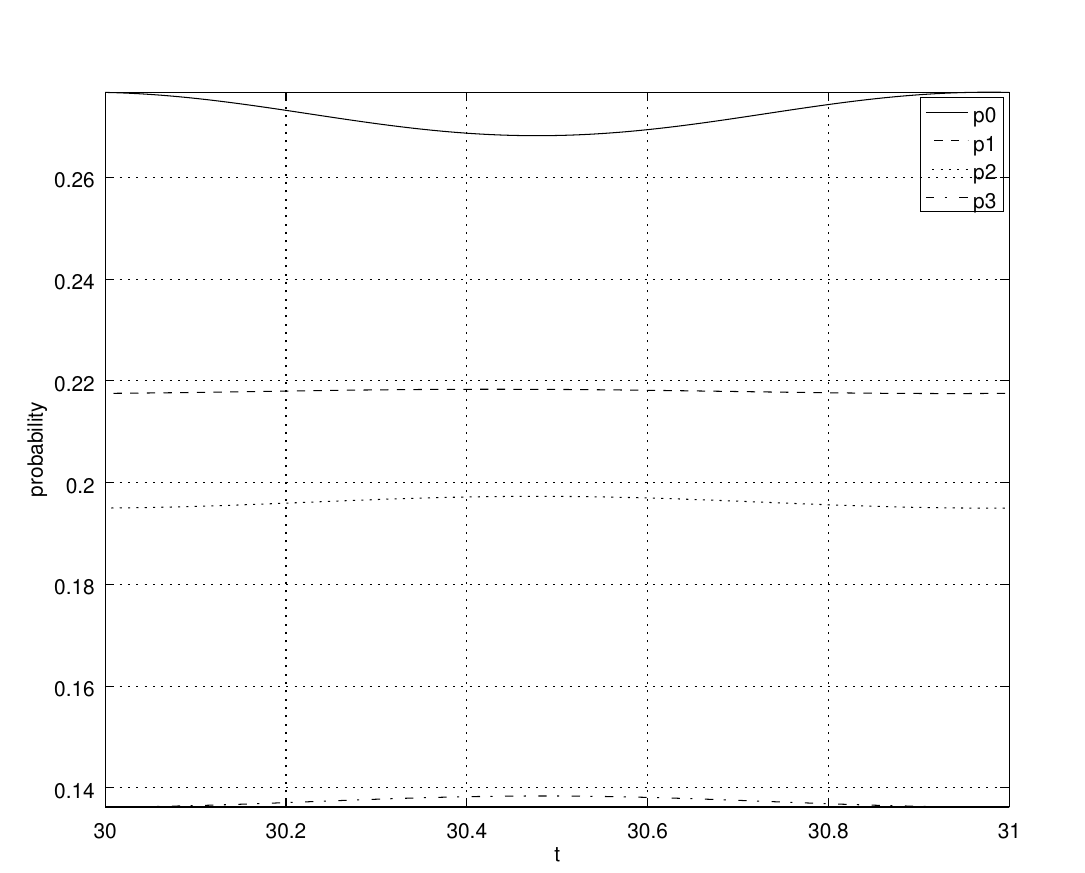}
%\vspace{-3cm}
\caption{The limiting probabilities $p_i(t)$, $0 \le i \le 3$, for $\mu =1$.}
%\Description{The plot for Example 2}
\label{fig:08}       % Give a unique label
\end{figure}

\begin{figure}[h]
\centering
  \includegraphics[width=0.5\linewidth]{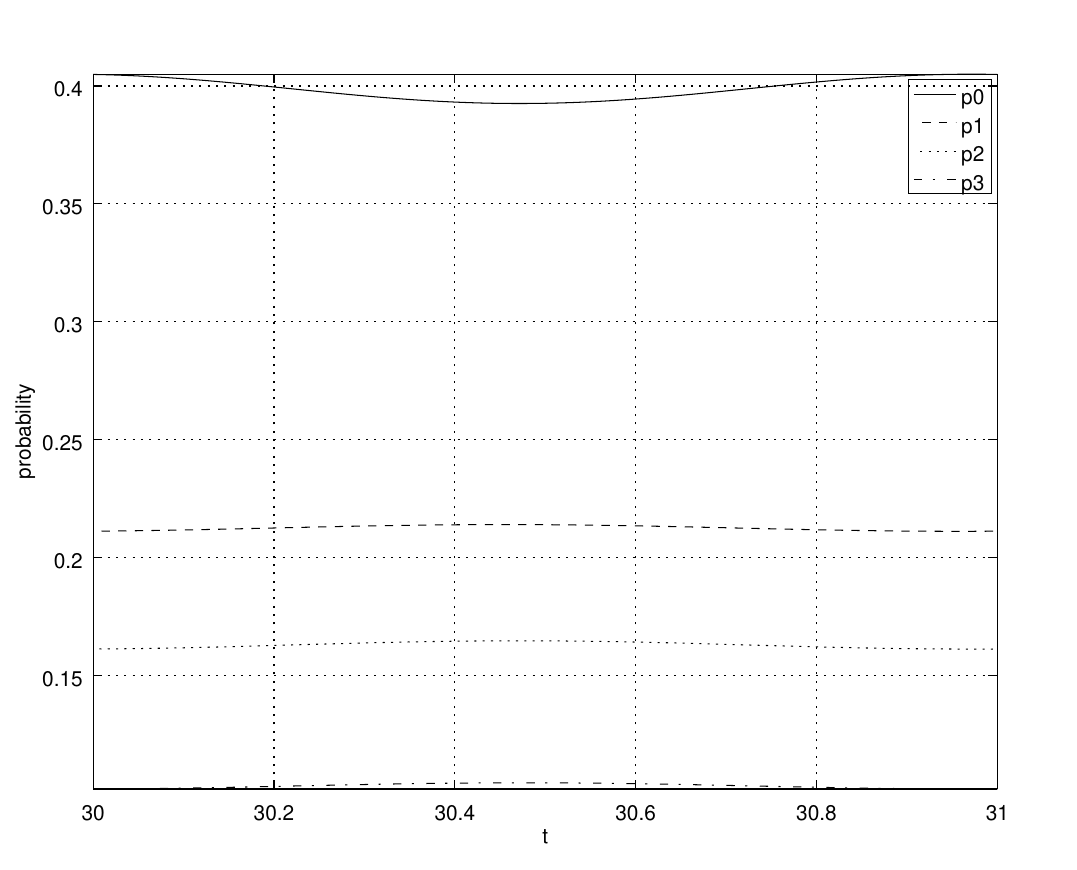}
%\vspace{-3cm}
\caption{The limiting probabilities $p_i(t)$, $0 \le i \le 3$, for $\mu =1.5$.}
%\Description{The plot for Example 2}
\label{fig:10}       % Give a unique label
\end{figure}

%%%%%%%%%%%%%%% 2

\begin{figure}[h]
\centering
  \includegraphics[width=0.5\linewidth]{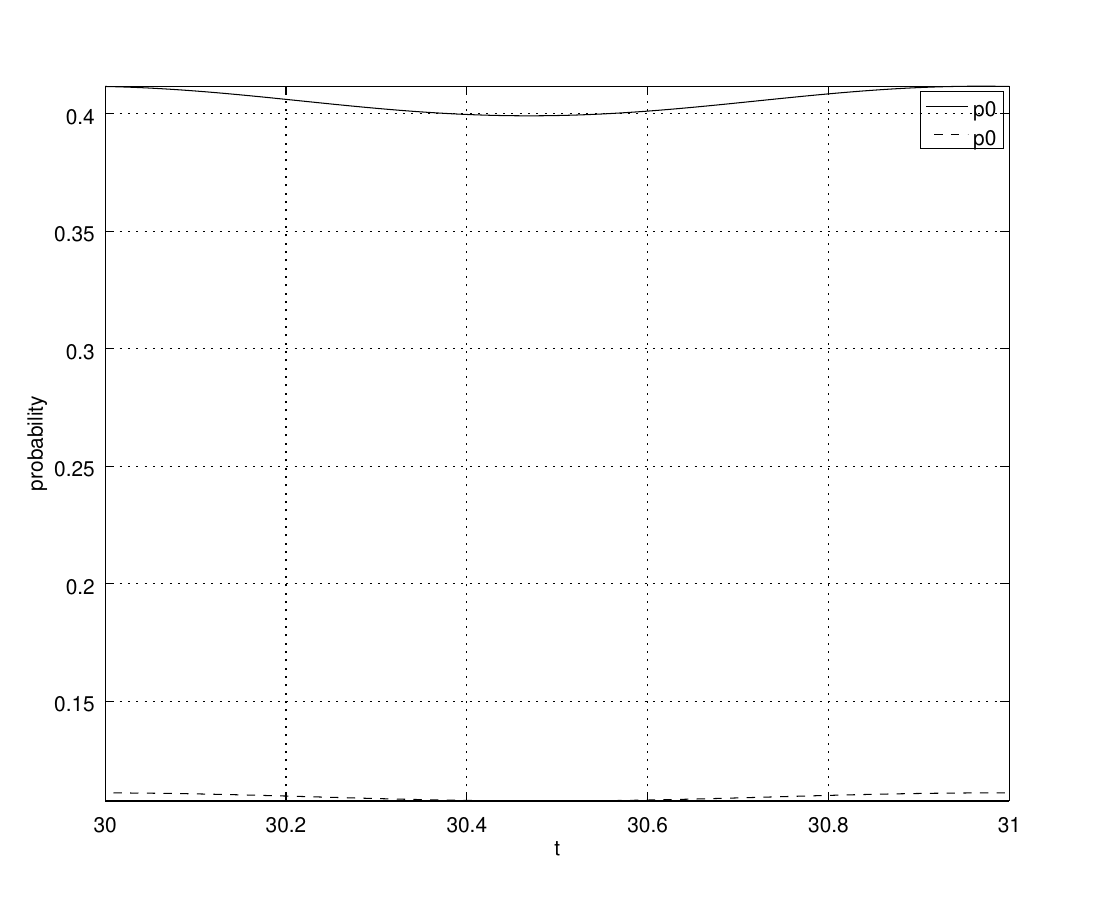}
%\vspace{-3cm}
\caption{The limiting probability $p_0(t)$  for $q=0.7$ (dotted line) and for $ q = 0.3 $ (solid line), $\lambda(t)=0.8+0.1 \sin (2 \pi t)$.}
%\Description{The plot for Example 3}
\label{fig:11}       % Give a unique label
\end{figure}

\begin{figure}[h]
\centering
  \includegraphics[width=0.5\linewidth]{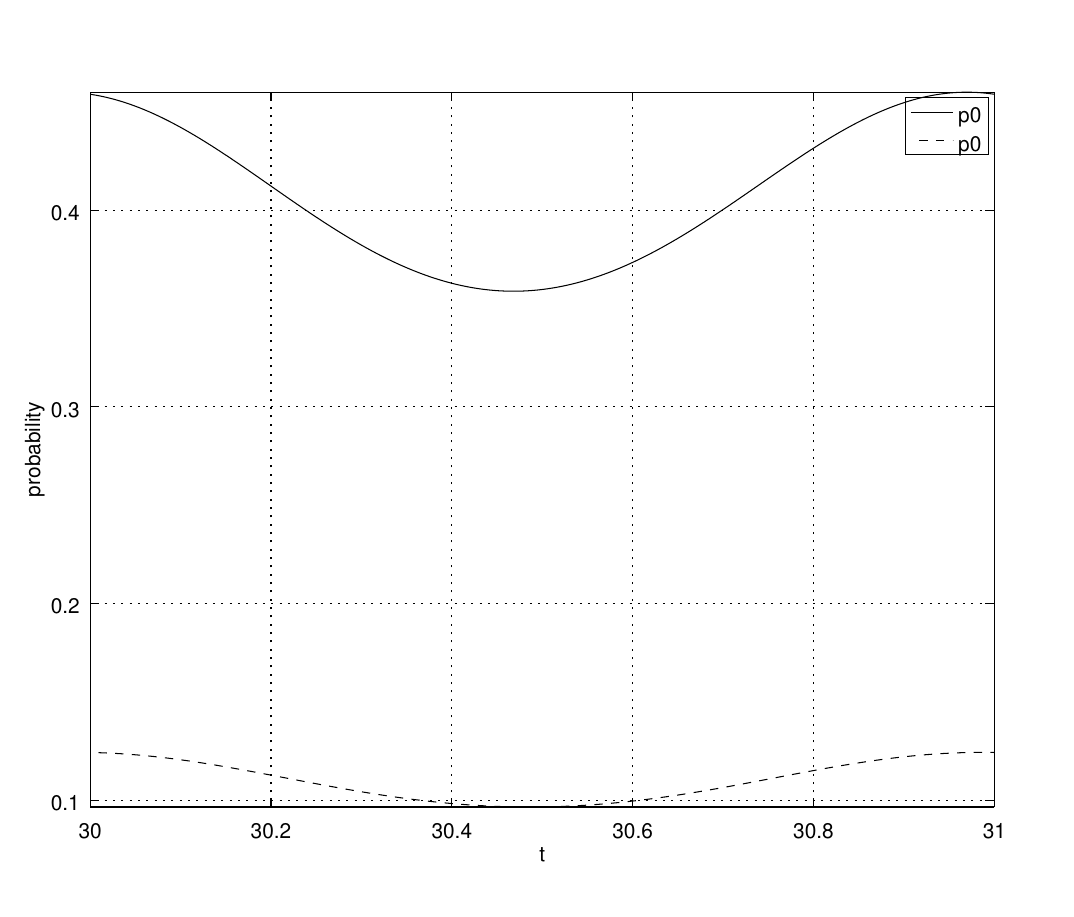}
%\vspace{-3cm}
\caption{The limiting probability $p_0(t)$ for $q=0.7$ (dotted line) and for $ q = 0.3 $ (solid line), $\lambda(t)=0.8+0.8 \sin (2 \pi t)$.}
%\Description{The plot for Example 3}
\label{fig:13}       % Give a unique label
\end{figure}

%%%%%%%%%%%%%%% 3

\begin{figure}[h]
\centering
  \includegraphics[width=0.5\linewidth]{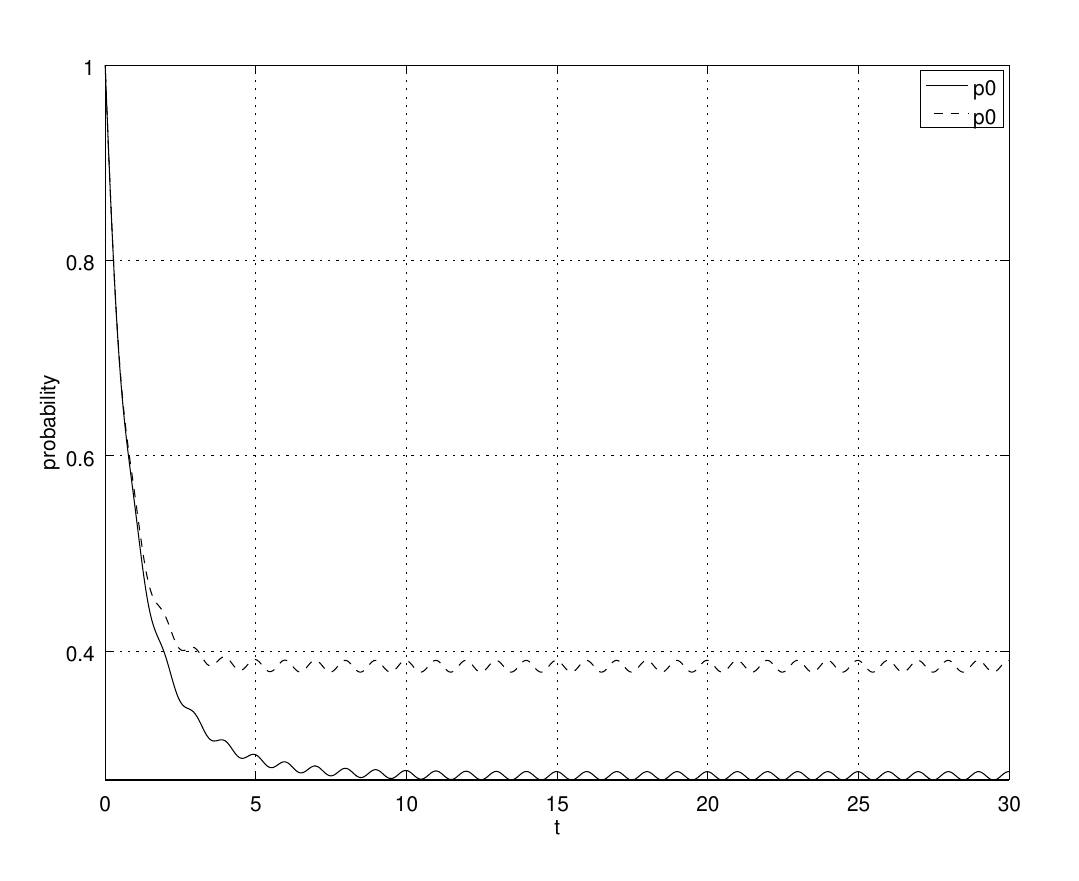}
%\vspace{-3cm}
\caption{Rate of convergence of the empty system
probability $p_0(t)$ for $q = 0.5$, $\mu =1$ (solid line) and for $q = 0$, $\mu =0.5$ (dotted line) in the interval $[0,30]$.}
%\Description{The plot for Example 4}
\label{fig:21}       % Give a unique label
\end{figure}

\begin{figure}[h]
\centering
  \includegraphics[width=0.5\linewidth]{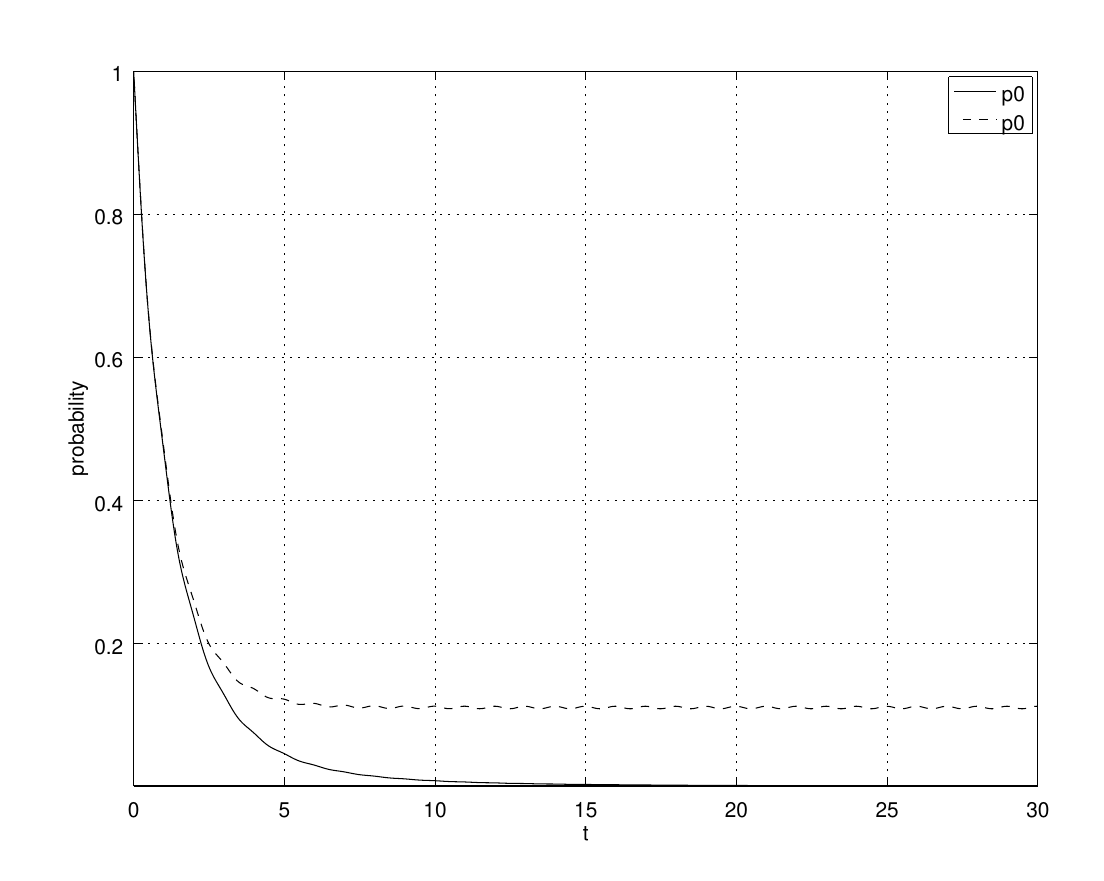}
%\vspace{-3cm}
\caption{Rate of convergence of the empty system
probability $p_0(t)$ for $q = 0.9$, $\mu =1$ (solid line) and for $q = 0$, $\mu =0.1$ (dotted line) in the interval $[0,30]$.}
%\Description{The plot for Example 4}
\label{fig:23}       % Give a unique label
\end{figure}

\end{document}